\documentclass[12pt]{article}
\usepackage{graphicx}
\usepackage{amsmath,bm,amssymb}
\usepackage{color}
\usepackage{amsthm}
\usepackage{authblk}

\def\hind{\hangindent=2pc\hangafter=1}
\newfont{\smcaps}{cmcsc10 scaled\magstep1}
\newcommand{\var}{{\rm ~Var\,}}

\newcommand{\E}{{\rm ~E\,}}

 at 11truept

\begin{document}
\baselineskip=22pt

\title{Improved Spread-Location Visualization}
\date{}
\author{A. Ian McLeod\\Western University}
\maketitle
\smallskip
\noindent McLeod, A.I. (1999),
Improved spread-location visualization.
{\it Journal of Graphical and Computational Statistics}  8/1, 135-141.\\
\hfill 

\noindent {\bf Abstract.}
Cleveland (1979, 1993) introduced the spread-location plot as a
diagnostic plot suitable for many types of fitted statistical models.
The spread-location plot
which plots the absolute residual or square-root absolute residual vs. fitted
value along with a robust loess smooth is
a useful replacement for the customary practice of
plotting residuals vs. fitted values.
In this note,
we show that neither absolute residual or square-root absolute residual is
always appropriate for error distributions likely to be encountered
in actual applications.
Instead we recommend the simple practice of examining the distribution of the
absolute residuals with a boxplot and choosing whatever transformation
is necessary to obtain symmetry before constructing the spread-location
plot.
We conclude with an illustrative example.

\newpage
\strut
\\{\bf 1. Introduction\hfill}
\bigskip
\break
It is standard practice in regression analysis to plot the residuals against
the fitted values.
If an increase in variability is detected in this plot
then a type of heteroscedasticity exists which can usually be removed
by a power transformation (Bartlett, 1947).
Cleveland (1979) pointed out that a drawback of this approach
is that simply more data at high fitted values can give the impression
of an increase in variability.
So Cleveland (1979) suggested plotting the absolute residuals.
In this plot an increase in variance corresponds to a monotonic trend.
Cleveland (1979) introduced the robust loess nonparametric regression smooth
to visually assess the trend in this plot.
Since the distribution of the absolute residuals are often skewed,
Cleveland (1993) recommends plotting the square root absolute
residual against fitted value.
Lack of skewness is important not only in improving our ability
to visually assess the spread-location plot
but it is also assumed when fitting a robust loess curve.
In this note we show that the square root transformation may not
always be the best choice and we recommend the simple expedient
of using a boxplot to choose a suitable transformation.

\strut
\\{\bf 2. Transformation to Symmetry of Absolute Residuals\hfill}
\bigskip
\break
Let $E$, the residual, be a random variable with mean zero and let $A=|E|$.
In many cases the distribution of $A$ will be positively skewed.
An obvious exception would be the case where $E$ has a uniform
distribution.
In general, we seek a power transformation of $A$,
$A^{(p)} = A^{p},\  p\ne 0$ and
$A^{(p)} = \log(A),\  p = 0$,
for which the lack of symmetry is reduced.
It is assumed that with probability 1, $A > 0$.

Consider the case where $E$ has probability density function $f(e)$
and we will assume that the third moment exists.
Then the density function of $A$ is $g(a) = f(a)+f(-a)$ which
reduces to $g(a) = 2 f(a)$ when $f$ is symmetric.
The density function for $B=A^{(p)}$ is given by
$h(b) = q b^{q-1} g(b)$, where $q = 1/p$ and $p \ne 0$.
When $p = 0$, $h(b) = e^b g(e^b)$.
Using a symbolic program such as {\it Mathematica\/},
we can easily evaluate the skewness coefficient of transformed variable $B$,
$$\sqrt \beta_1 = { \E\{ (B- E\{B\})^3 \} \over \var(B)^{3/2}},$$
for any particular distribution of $E$.
While only good numerical integration techniques are required to get
the results shown in Table 1, the computations are more easily performed in
a symbolic computation environment where sophisticated symbolic integration
routines are available.
Most of the results in Table 1 were obtained symbolically.
For convenience the decimal approximations are given.
A {\it Mathematica} notebook containing the detailed derivations
of all results in Table 1 is available from my homepage.

\begin{figure}
\begin{center}
{\normalsize Table 1. Skewness Coefficient of $|E|^p$} \\
\bigskip
\begin{tabular}{lrrrrrr} \hline  \hline
\noalign{\vskip6pt}
  Distribution of $E$
& $p=1$
& $p=0.5$
& $p=0.4$
& $p=0.3\dot 3$
& $p=0.25$
& $p=0$
\\
\noalign{\vskip6pt}
\hline
\noalign{\vskip6pt}
\strut
Gaussian &$0.995$ &$0.0841$ &$-0.1344$ &$-0.2948$ &$-0.5192$ &$-1.535$\\
$t_5$ &$2.550$ &$0.5827$ &$0.2766$ &$0.0687$ &$-0.2050$ &$-1.313$\\
5\% 3 $\sigma$ contamination &$2.674$ &$0.5872$ &$0.2520$ &$0.0286$ &$-0.2600$ &$-1.140$\\
Laplace &$2.000$ &$0.6311$ &$0.3586$ &$0.1681$ &$-0.0872$ &$-1.395$\\
\noalign{\vskip6pt}
\hline
\end{tabular}
\end{center}
\end{figure}
As shown in Table 1,
for the normal distribution, the square root transformation
suggested by Cleveland (1993) is appropriate.
But for other distributions, there are better choices.
For example, fat-tailed distributions often occur in practice
(Tukey, 1960) and Table 1 shows that
if the errors are approximated either by a $t$-distribution on
5 df or by a contaminated normal distribution with a
5\% probability of a scale inflation by a factor of 3,
then a cube root transformation does better.
For Laplace errors, $p=0.25$ does better.
\strut
\\{\bf 3. Illustrative Example\hfill}
\medskip
\break
The ethanol data of Cleveland (1993, p.217, Figure 4.33)
provides as a nice illustration.
The boxplots in Figure 1 below show that a log transformation, $p = 0$,
makes the absolute error distribution approximately symmetric.
Figure 2 compares the spread-location plots using a square-root
and log transformation.
Notice that the plot with the log transformation is more symmetrically distributed about the
fitted curve.
Since the amount of variation exhibited in the data is still
quite large, compared to the curve, it is safe to assume
that monotone spread is not present.

In practice, if a large degree of monotone spread is detected,
it may be desirable to use a slightly different power transformation
than that initially chosen
on the absolute residuals in order to get symmetric deviations
from the loess trend on the spread-location visualization.
\medskip
\begin{center}
Figures 1 and 2 here
\end{center}

\newpage
\begin{center}
{\bf REFERENCES\hfill}
\end{center}
\medskip
\parindent 0pt

\hind
Bartlett, M.S. (1947),
``The Use of Transformations'',
{\it Biometrics\/}, 3, 37--52.

\hind
Cleveland, W.~S. (1979),
``Robust Locally Weighted Regression and Smoothing Scatterplots'',
{\it Journal of the American Statistical Association\/},
74, 829--836.

\hind
Cleveland, W.~S. (1993),
{\it Visualizing Data\/}. Summit, New Jersey: Hobart Press.

\hind
Tukey, J.\ W. (1960),
``A survey of sampling from contaminated distributions''
in Contributions to Probability and Statistics: Essays in Honor of
Harold Hotelling, Edited by I. Olkin, S.\ G. Ghurye, W. Hoeffding,
W.\ G. Madow and H.\ B. Mann, Standford University Press, Standford.

\strut
\vspace{5in}
\vfill\eject
\begin{figure}
\caption{Boxplots absolute residuals and transformed absolute residuals.
Data rescaled to fit along common axis.}
\begin{center}
\includegraphics[width=0.8\textwidth]{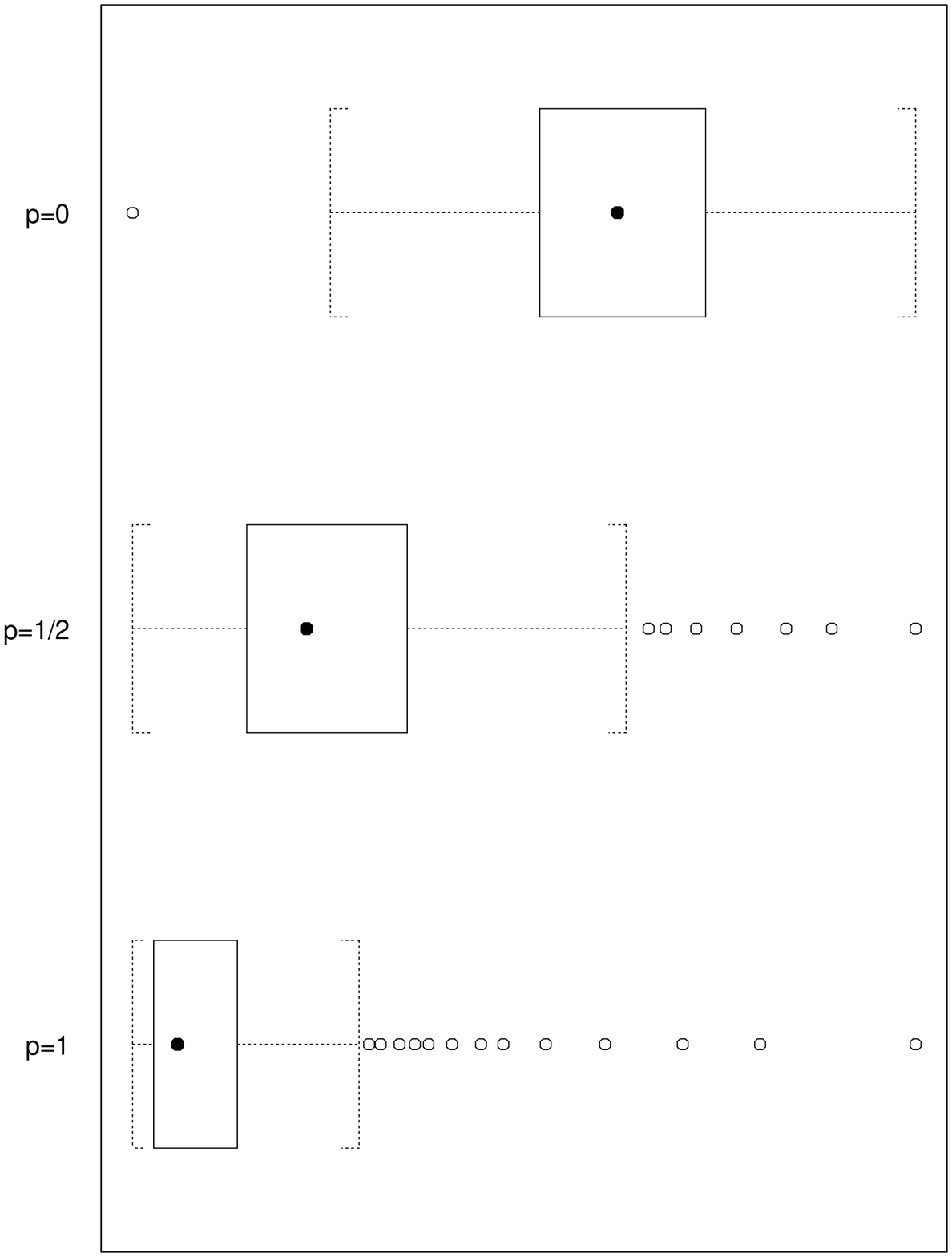}
\end{center}
\label{fig:boxplots}
\end{figure}


\newpage

\begin{figure}
\caption{Boxplots absolute residuals and transformed absolute residuals.
Data rescaled to fit along common axis.}
\begin{center}
\includegraphics[width=0.8\textwidth]{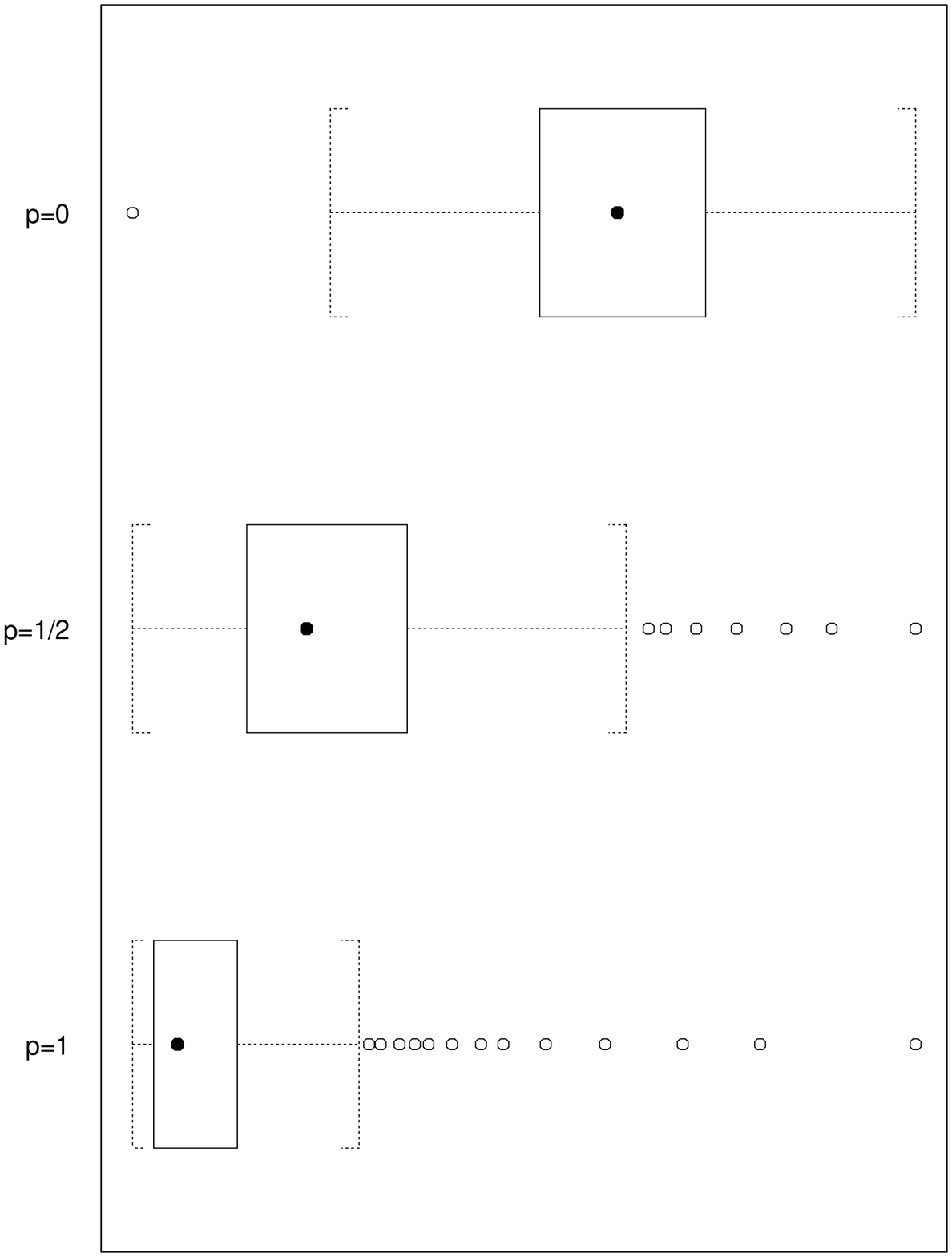}
\end{center}
\label{fig:ethanol}
\end{figure}


\end{document}